\input amstex
\documentstyle{amsppt}
\magnification=\magstep1                        
\hsize6.5truein\vsize8.9truein                  
\NoRunningHeads
\loadeusm

\magnification=\magstep1                        
\hsize6.5truein\vsize8.9truein                  
\NoRunningHeads
\loadeusm

\document
\topmatter

\date February 6, 2016 
\enddate

\title Inequalities for exponential sums
\endtitle

\rightheadtext{Inequalities for exponential sums}

\author Tam\'as Erd\'elyi
\endauthor

\address Department of Mathematics, Texas A\&M University, College Station, Texas 77843 \endaddress

\email terdelyi\@math.tamu.edu (T. Erd\'elyi)
\endemail

\email terdelyi\@math.tamu.edu
\endemail

\endthanks


\abstract
We study the classes
$${\Cal E}_n := \left \{f: f(t) = \sum_{j=1}^n{a_j e^{\lambda_jt}}\,, \enskip \enskip a_j, \lambda_j \in {\Bbb C} \right \}\,,$$
$${\Cal E}_n^+ := \left \{f: f(t) = \sum_{j=1}^n{a_j e^{\lambda_jt}}\,, \enskip \enskip a_j, \lambda_j \in {\Bbb C}\,, \enskip \text {\rm Re}(\lambda_j) \geq 0 \right \}\,,$$
and
$${\Cal T}_n := \left \{f: f(t) = \sum_{j=1}^n{a_j e^{i\lambda_jt}}\,, \enskip a_j \in {\Bbb C}, \enskip \lambda_1 < \lambda_2 < \cdots < \lambda_n \right \}\,.$$
A highlight of this paper is the asymptotically sharp inequality 
$$|f(0)| \leq (1 + \varepsilon_n)\, 3n \, \|f(t)e^{-9nt/2}\|_{L_2[0,1]}, \qquad f \in {\Cal T}_n \,,$$
where $\varepsilon_n$ converges to $0$ rapidly as $n$ tends to $\infty$. 
The inequality  
$$\sup_{f \in {\Cal T}_n}{ \frac{|f(0)|}{\|f\|_{L_2{[0,1]}}}} \geq n\,.$$
is also observed. Our results improve an old result of G. Hal\'asz and a recent result of G. K\'os. 
We prove several other essentially sharp related results in this paper.
\endabstract

\endtopmatter

\head 1. Introduction and Notation \endhead
The well known results of Nikolskii assert that the essentially sharp inequality
$$\|P\|_{L_q[-1,1]} \leq c(p,q)n^{2/p-2/q}\|P\|_{L_p[-1,1]}$$
holds for all algebraic polynomials $P$ of degree at most $n$
with complex coefficients and for all $0 < p \leq q \leq \infty$, while the essentially sharp inequality
$$\|Q\|_{L_q[-\pi,\pi]} \leq c(p,q)n^{1/p-1/q}\|Q\|_{L_p[-\pi,\pi]}$$
holds for all trigonometric polynomials $Q$ of degree at most $n$
with complex coefficients and for all $0 < p \leq q \leq \infty$. The subject
started with two remarkable papers, [24] and [28]. There are quite a few related papers
in the literature, and several books discuss inequalities of this
variety with elegant proofs. See [3] and [12], for example.
In this paper we focus on the classes
$${\Cal E}_n := \left \{f: f(t) = \sum_{j=1}^n{a_j e^{\lambda_jt}}\,, \enskip \enskip a_j, \lambda_j \in {\Bbb C} \right \}\,,$$
$${\Cal E}_n^+ := \left \{f: f(t) = \sum_{j=1}^n{a_j e^{\lambda_jt}}\,, \enskip \enskip a_j, \lambda_j \in {\Bbb C}\,, \enskip \text {\rm Re}(\lambda_j) \geq 0 \right \}\,,$$
$${\Cal E}_n^- := \left \{f: f(t) = \sum_{j=1}^n{a_j e^{\lambda_jt}}\,, \enskip \enskip a_j, \lambda_j \in {\Bbb C}\,, \enskip \text {\rm Re}(\lambda_j) \leq 0 \right \}\,,$$
and
$${\Cal T}_n := \left \{f: f(t) = \sum_{j=1}^n{a_j e^{i\lambda_jt}}\,, \enskip a_j \in {\Bbb C}\,, \enskip \lambda_1 < \lambda_2 < \cdots < \lambda_n \right \}\,.$$
These classes were studied in several publications. See [20], [22], [23], and [29], for example.
For the sake of brevity let
$$\|f\|_A := \sup_{t \in A}{|f(t)|}$$
for a complex-valued function $f$ defined on a set $A \subset {\Bbb R}$.
Section 19.4 of Tur\'an's book [29] refers to the following result of G. Hal\'asz: 
$$|f(0)| \leq cn^5 \|f\|_{[0,1]}\,, \qquad f \in {\Cal E}_n^+\,,$$
where $c>0$ is an absolute constant. This was improved recently by G. K\'os [20] to 
$$|f(0)| \leq 10 \, \frac{5n}{5n-1} n^2 \|f\|_{L_1[0,1]}\,, \qquad f \in {\Cal E}_n^+\,, \tag 1.1$$
where $cn^2$ is the best possible size of the factor in this inequality.
He also proved that
$$|f(0)| \leq 2n \|f\|_{L_2[0,1]}\,, \qquad f \in {\Cal E}_n^+\,, \tag 1.2$$
where $cn$ is the best possible size of the factor in this inequality.
The technique used in [20] is based on integrating discrete inequalities similar 
to Tur\'an's first and second main theorems in the theory of power sums. This technique 
was also used by Tijdeman as it was demonstrated, for example in Section 27 of Tur\'an's book [29].
This answers a question of S. Denisov asked from me in e-mail communications. I was not aware 
of the above results when I started to write this paper. In this paper we recapture the above 
inequalities with better constants for all $f \in {\Cal T}_n$. Namely we prove that
$$|f(0)| \leq cn^2 \|f\|_{L_1[0,1]}\,, \qquad f \in {\Cal T}_n\,, \tag 1.3$$
with $c = 2 + \log 4 + \varepsilon_n = 3.3862\ldots$ and
$$|f(0)| \leq cn \|f\|_{L_2[0,1]}\,, \qquad f \in {\Cal T}_n\,, \tag 1.4$$
with $c = (2 + \log 4 + \varepsilon_n)^{1/2} = 1.8401\ldots$.
S. Denisov [11] has just proved that the constant $c = (2 + \log 4 + \varepsilon_n)^{1/2}  = 1.8401\ldots$ 
can be further improved to $c = \pi/2 = 1.5707\ldots$ in (1.4). 
Denisov's approach also uses a Hal\'asz-like construction first, 
which may be found in [19] and it also appears as Lemma 10.8 in [29], but after that it employs 
a duality argument and an old result of Lachance, Saff, and Varga [21], which is not used by K\'os.   
We note that Denisov's improvement of (1.4) can also be seen for all $f \in {\Cal E}_n^+$ by 
modifying K\'os's approach. Indeed, it is proved in [21] that 
$$\sigma_k := \min \big\{\|P(e^{it})\|_{[0,2\pi]}: P(0) = 1, \enskip P(1) = 0, \enskip P \in {\Cal P}_k^c \big\} = 
\left( \sec \frac{\pi}{2(k+1)} \right)^{k+1}\,,$$ 
where ${\Cal P}_k^c$ denotes the set of all algebraic polynomials of degree at most $k$ with complex 
coefficients. Hence there are polynomials $H_k \in {\Cal P}_k^c$ such that $H_k(0) = 1$ and 
$$\|H_k(e^{it})\|_{[0,2\pi]} \leq  \left( \sec \frac{\pi}{2(k+1)} \right)^{k+1} = \exp \left( \frac{\pi^2}{8k} + O\left( \frac{1}{k^2} \right) \right)\,.$$ 
Using the above $H_k \in {\Cal P}_k^c$ instead of the $H_k \in {\Cal P}_k^c$ in K\'os's proof satisfying only 
$$\|H_k(e^{it})\|_{[0,2\pi]} \leq \exp \left( \frac 2k \right) \,,$$
we get Denisov's improvement of (1.4) can be extended to all $f \in {\Cal E}_n^+$, that is, 
$$|f(0)| \leq \frac{\pi n}{2} \, \|f\|_{L_2[0,1]}\,, \qquad f \in {\Cal E}_n^+\,. \tag 1.5$$
In Section 3 the infinite-finite range inequality is stated  
$$\int_0^{\infty}{|f(t)|^2 e^{-t} \, dt} \leq (1+\varepsilon_n)^2 \, \int_0^{9n}{|f(t)|^2 e^{-t} \, dt}$$
for every $f \in {\Cal E}_n^-$, in particular, for every $f \in {\Cal E}_n^-$.
Here $(1 + \varepsilon_n)^2 := 1 + 8190e^{-n/10}$ is an appropriate choice.
As a consequence we prove that 
$$|f(0)| \leq (1 + \varepsilon_n) \, 3n \, \|f(t)e^{-9t/2}\|_{L_2[0,1]}\,, \qquad f \in {\Cal T}_n\,,$$
where $\varepsilon_n$ is the same as before, 
and for every $\lambda_1 < \lambda_2 < \cdots < \lambda_n$ there is an $f \in {\Cal T}_n$ of the form 
$$f(t) = \sum_{j=1}^n{a_j e^{i\lambda_jt}}\,, \qquad a_j \in {\Bbb C}\,, \tag 1.6$$
such that
$$|f(0)| >  3n \, \|f(t)e^{-9nt/2}\|_{L_2[0,1]}\,.$$
Other Nikolskii-type inequalities comparing the $L_p[0,1]$ and $L_q[0,1]$ norms of exponential 
sums $f \in {\Cal T}_n$ are also established in Section 2. We use quite different techniques based 
on the knowledge of M\"untz-Legendre orthonormal polynomials studied in [8] and Section 3.4 of [3]. 
We obtain interesting Markov-type inequalities as well for the derivatives of exponential sums 
$f \in {\Cal T}_n$, but such a Markov-type inequality cannot depend only on $n$, it depends  on 
the exponents $\lambda_1 < \lambda_2 < \cdots < \lambda_n$. We also examine how far our 
estimates are from being sharp, and it turns out that our main results proved in this paper 
are essentially sharp. Most importantly, the inequality 
$$\sup_{f \in {\Cal T}_n}{ \frac{|f(0)|}{\|f\|_{L_2{[0,1]}}}} \geq n\,.$$
is also observed in Section 2. The inequality   
$$|f(0)| \leq n\|f\|_{L_2[0,1]}$$ 
for every $f \in {\Cal E}_n^+$ of the form
$$f(t) = \sum_{j=1}^n{a_j e^{\lambda_j t}}\,, 
\qquad a_j \in {\Bbb R}\,, \enskip 0 \leq \lambda_1 < \lambda_2 < \cdots < \lambda_n\,, \tag 1.7$$
is stated in Section 4. This inequality is sharp. We suspect that the above inequality holds for all 
$f \in {\Cal T}_n$ or perhaps for all $f \in {\Cal E}_n^+$ at least with $n$ 
replaced by $(1+\varepsilon_n)n$, where $\varepsilon_n$ tends to $0$ as $n$ tends to $\infty$. 
Markov-Nikolskii-type inequalities for ${\Cal T}_n$ are established in Section 5. 
Markov-Nikolskii-type inequalities for $f \in {\Cal E}_n$ with nonnegative exponents are formulated  
in Section 6. We claim that
$$|f^{\prime}(0)| \leq (1 + \varepsilon_n) \, 3^{-1/2} \, n^3 \|f\|_{L_2[0,1]}$$
for every $f \in {\Cal E}_n^+$ of the form (1.7),  
where the quantity $\varepsilon_n$ (determined exactly in the proof) tends to $0$ an $n$ 
tends to $\infty$. This inequality is sharp. 
Section 7 offers an essentially sharp pointwise Nikolskii-type inequality for ${\Cal E}_n$, 
namely we claim that
$$\left( \frac{(n-2)\log 2}{4\min \{y-a,b-y\}} \right)^{1/2} \leq \sup_{f \in {\Cal E}_n}{\frac{|f(y)|}{\|f\|_{L_2[a,b]}}} 
\leq \left( \frac{2n}{\min\{{y-a,b-y\}}} \right)^{1/2}$$
for every $y \in (a,b)$. In Section 8 we offer the Bernstein-type inequality
$$|f^{\prime}(0)| \leq 2e(\lambda + n + 1) \, \|f\|_{[-1,1]}$$
for every $f \in {\Cal T}_n$ of the form (1.6), where
$$\lambda := \max_{1 \leq j \leq n} {|\lambda_j|}\,, 
\qquad \lambda_1 < \lambda_2 < \cdots < \lambda_n\,. \tag 1.8$$
This inequality is sharp up to the factor $2e$. Namely, for every real number $\lambda > 0$ 
and integer $n \geq 1$ there is an $f \in {\Cal T}_n$ of the form (1.6) with (1.8) such that
$$|f^{\prime}(0)| \geq \frac 14 \, (\lambda + n - 3) \, \|f\|_{[-1,1]}\,.$$
In Section 9 the Markov-type inequality 
$$\|f^{\prime}\|_{[0,1]} \leq (1 +\varepsilon_n) 
\left( 108n^5 + \sum_{k=1}^n{\lambda_k^2} \right)^{1/2} \|f\|_{[0,1]}$$
for every $f \in {\Cal T}_n$ of the form (1.6) is established, where the quantity $\varepsilon_n$ (
determined exactly in the proof) tends to $0$ an $n$ tends to $\infty$. We record an observation 
showing how far the above Markov inequality is from being sharp. 
Markov-type inequalities for ${\Cal E}_n^-$ and ${\Cal T}_n$ in $L_2[0,\infty)$ with the
Laguerre weight are established in Section 10. Our Theorem 10.1 extends Lubinsky's Theorem 3.2 
in [22] from the case of exponential sums with purely imaginary exponents to the case of exponential 
sums with complex exponents. Our only result in Section 11 is a version of Theorem 10.1, a Markov-type 
inequality for ${\Cal E}_n^-$ in $L_2[0,\infty)$ without a weight. We prove our new results in Section 13. 
Lemmas needed in the proofs of our new results are stated and proved in Section 12. Combining Tur\'an's 
power sum method with results in [9], [10], and [17], we may be able to prove other interesting results 
in the future. We close the paper with an Appendix listing results closely related to our new results in 
this paper. Theorems 14.1--14.6 have been proved by subtle Descartes system methods which can be employed 
only in the case of exponential sums with real exponents but not in the case of complex exponents. The 
reader may find it useful to compare the results in Section 14 with the new results of the paper.

\medskip

\head 2. New Results: Nikolskii-type inequalities for ${\Cal T}_n$ \endhead

Answering a question by Sergey Denisov (e-mail communications) in this paper we prove the following
new results. Observe that while our constant $(8+\varepsilon_n)$ is not as good as $\pi/2$ or even $2$, 
there is a rapidly decreasing weight function $w(t) = e^{-nt}$ pushing down the $L_2[0,1]$ norm at the 
right-hand side. 

\proclaim{Theorem 2.1} 
We have
$$|f(0)| \leq (8+\varepsilon_n)^{1/2} \, n \, \|f(t)e^{-nt}\|_{L_2[0,1]}, \qquad f \in {\Cal T}_n\,,$$
where $(8 + \varepsilon_n)^{1/2} := 8^{1/2}(1+2e^{-2n})^{1/2}$, 
and for every $\lambda_1 < \lambda_2 < \cdots < \lambda_n$ there is an $f \in {\Cal T}_n$ of the form (1.6) such that
$$|f(0)| > 8^{1/2}  \, n \, \|f(t)e^{-4nt}\|_{L_2[0,1]}\,.$$
\endproclaim

Our next theorem recaptures K\'os's inequality (1.1) with a 
a constant better than $c=2$ but not as good as $c=\pi/2$. The constant specified in our theorem below 
seems to be the limit of what our essentially different method based on the explicit form of 
M\"untz-Legendre orthonormal polynomials gives.   

\proclaim{Theorem 2.2}
Let $\gamma_0 := 2 + \log 4 < \gamma \leq 4$. We have 
$$|f(0)| \leq (\gamma + \varepsilon_n)^{1/2} \, n \, \|f\|_{L_2[0,1]}, \qquad f \in {\Cal T}_n\,,$$
where 
$$(\gamma + \varepsilon_n)^{1/2} = \gamma^{1/2}(1+\delta^{-2} e^{-\delta \gamma n})^{1/2}\,, \qquad 
\delta := \frac{\gamma - \gamma_0}{8}\,.$$
\endproclaim

Observe that if $f \in {\Cal T}_n$ and $g(t) = f(-t)$, then $g \in {\Cal T}_n$. Hence the 
extension of Theorem 2.2 formulated by our next couple of theorems follows easily.   
 
\proclaim{Theorem 2.3}
We have
$$\|f\|_{[0,1]} \leq \frac{\pi n}{2} \, \|f\|_{L_2[0,1]}\,, \qquad f \in {\Cal T}_n\,.$$
\endproclaim

\proclaim{Theorem 2.4}
We have 
$$\|f\|_{[0,1]} \leq  \, \left( \frac{\pi n}{2} \right)^{2/q} \|f\|_{L_q[0,1]}\,, \qquad f \in {\Cal T}_n\,, \quad q \in (0,2]\,.$$
\endproclaim

\proclaim{Theorem 2.5}
We have
$$\|f\|_{L_p[0,1]} \leq \left( \frac{\pi n}{2} \right)^{2/q - 2/p} \|f\|_{L_q[0,1]}\,, \qquad f \in {\Cal T}_n\,, \quad 0 < q < p \leq \infty\,, \quad q \leq 2\,.$$
\endproclaim

Note that the case $q=1$ of Theorem 2.4 improves K\'os's inequality (1.1) to 
$$\|f\|_{[0,1]} \leq \frac {\pi^2n^2}{4} \, \|f\|_{L_1[0,1]}\,, \qquad f \in {\Cal T}_n\,.$$ 

\proclaim{Theorem 2.6} We have
$$\sup_{f \in {\Cal T}_n}{ \frac{|f(0)|}{\|f\|_{L_2{[0,1]}}}} \geq n\,.$$
\endproclaim

\proclaim{Theorem 2.7} There is an absolute constant $c > 0$ such that
$$\sup_{f \in {\Cal T}_n}{ \frac{|f(0)|}{\|f\|_{L_q{[0,1]}}}} \geq c^{1+1/q}(1+qn)^{2/q}, \qquad q \in (0,\infty)\,.$$
\endproclaim

\noindent {\bf Remark 2.8.}
It remains open what are the right extensions of Theorems 2.4 and 2.5 to $q > 2$. Note that Theorem 1.8 implies that
$$\sup_{f \in {\Cal T}_n}{ \frac{|f(0)|}{\|f\|_{L_q{[0,1]}}}} \geq \sup_{f \in {\Cal T}_n}{ \frac{|f(0)|}{\|f\|_{L_q{[-\pi,\pi]}}}} \geq c_q n^{1/2}, \qquad q \in (0,\infty)\,,$$
with a constant $c_q > 0$ depending only on $q > 0$. Hence the right upper bound in Theorem 2.4 is 
somewhere between $c_q n^{1/2}\,$ and $(\gamma + \varepsilon)^{1/2} \, n$  for all $q > 2$. In particular, 
Theorems 2.4 cannot remain true for $q > 4$.

\medskip

\proclaim{Theorem 2.9}
We have
$$|f(0)| \leq (8+\varepsilon_n)^{1/2} \, c_q n^{1/2+1/q}  \|f(t)e^{-nt}\|_{L_q[0,1]}\,, \quad f \in {\Cal T}_n\,, \quad q \in (2,\infty)\,,$$
where $\varepsilon_n$ is the same as in Theorem 2.1 and 
$$c_q := \left( \frac{q-2}{2q} \right)^{(q-2)/(2q)}\,.$$ 
\endproclaim

\proclaim{Theorem 2.10}
We have
$$\|f\|_{[0,1]} \leq (8+\varepsilon_n)^{1/2} c_q n^{1/2+1/q} \|f\|_{L_q[0,1]}\,, \quad f \in {\Cal T}_n\,, \quad q \in (2,\infty)\,,$$
where $\varepsilon_n$ is the same as in Theorem 2.1 and $c_q$ is the same as in Theorem 2.9. 
\endproclaim

\medskip

\head 3. New Results: An infinite-finite range inequality for ${\Cal E}_n^-$ with an application \endhead

Our next theorem is an infinite-finite range inequality for all $f \in {\Cal E}_n^-$. 

\proclaim{Theorem 3.1}
We have
$$\int_0^{\infty}{|f(t)|^2 e^{-t} \, dt} \leq (1+\varepsilon_n)^2 \, \int_0^{9n}{|f(t)|^2 e^{-t} \, dt}$$
for every $f \in {\Cal E}_n^-$, in particular, for every $f \in {\Cal T}_n\,.$ 
Here $(1 + \varepsilon_n)^2 := 1 + 8190e^{-n/10}$ is an appropriate choice.
\endproclaim

The theorem below establishes an asymptotically sharp version of K\'os's inequality 
$|f(0)| \leq 2n \|f\|_{[0,1]}$ in the presence of the rapidly decreasing weight function 
$w(t) = e^{-9nt/2}$ pushing down the $L_2[0,1]$ norm at the right-hand side. 

\proclaim{Theorem 3.2}
Let $\varepsilon_n$ be the same as in Theorem 2.11. We have
$$|f(0)| \leq (1 + \varepsilon_n) \, 3n \, \|f(t)e^{-9nt/2}\|_{L_2[0,1]}\,, \qquad f \in {\Cal T}_n\,,$$ 
and for every $\lambda_1 < \lambda_2 < \cdots < \lambda_n$ there is an $f \in {\Cal T}_n$ of the form 
(1.6) such that
$$|f(0)| >  3n \, \|f(t)e^{-9nt/2}\|_{L_2[0,1]}\,.$$
\endproclaim

\medskip

\head 4. A sharp Nikolskii-type inequality for $f \in {\Cal E}_n$ with nonnegative exponents \endhead 

Our next theorem establishes the best constant in the inequality $|f(0)| \leq cn \|f\|_{L_2[0,1]}$
for functions $f$ in a subclass ${\Cal E}_n$.

\proclaim{Theorem 4.1}
We have
$$|f(0)| \leq n\|f\|_{L_2[0,1]}$$
for every $f \in {\Cal E}_n^+$ of the form 
$$f(t) = \sum_{j=1}^n{a_j e^{\lambda_j t}}\,, \qquad a_j \in {\Bbb R}\,, \enskip 
0 \leq \lambda_1 < \lambda_2 < \cdots < \lambda_n\,.$$ 
This inequality is sharp. 
\endproclaim

\medskip

\head 5. New Results: Markov-Nikolskii-type inequalities for ${\Cal T}_n$ \endhead

The next theorem establishes the right result when $|f(0)|$ is replaced by $|f^{\prime}(0)|$ 
in Theorem 3.2.

\proclaim{Theorem 5.1}
Let $\varepsilon_n$ be the same as in Theorem 3.1. We have
$$|f^{\prime}(0)| \leq 27 \, (1 +\varepsilon_n) \, n^{3/2} 
\left( \sum_{k=1}^n{\left( \left( \frac{\lambda_k}{9n} \right)^2 + (k-1)^2 \right)} \right)^{1/2} 
\|f(t)e^{-9nt/2}\|_{L_2[0,1]}$$
for every $f \in {\Cal T}_n$ of the form (1.6), 
and for every $\lambda_1 < \lambda_2 < \cdots < \lambda_n$ there is an $f \in {\Cal T}_n$ of the form (1.6) such that
$$|f^{\prime}(0)| > 27\,n^{3/2} 
\left( \sum_{k=1}^n{\left( \left( \frac{\lambda_k}{9n} \right)^2 + (k-1)^2 \right)} \right)^{1/2}
\|f(t)e^{-9nt/2}\|_{L_2[0,1]}\,.$$
\endproclaim

The next theorem establishes the right result when $|f^{\prime}(0)|$ is replaced by 
$\|f^{\prime}\|_{[0,1]}$ in Theorem 5.1.

\proclaim{Theorem 5.2}
Let $\varepsilon_n$ be the same as in Theorem 3.1. We have
$$\|f^{\prime}\|_{[0,1]} \leq 27 \, (1 +\varepsilon_n) \, n^{3/2} 
\left( \sum_{k=1}^n{\left( 2\left( \frac{\lambda_k}{9n} \right)^2 + 8(k-1)^2 \right)} \right)^{1/2}
\|f\|_{L_2[0,1]}$$
for every $f \in {\Cal T}_n$ of the form (1.6). 
\endproclaim

To formulate our next observation, given $n \in {\Bbb N}$ and $\eta > 0$, we introduce 
the classes 
$${\Cal T}_n(\eta) := \Biggl\{f : f(t) = \sum^n_{j=1} {a_j e^{i\lambda_j t}}\,, \enskip a_j \in {\Bbb C}, \quad 0 < \lambda_1 < \lambda_2 < \cdots < \lambda_n < \eta \Biggr\}\,.$$

\proclaim{Theorem 5.3} We have 
$$\sup_{0 \neq f \in {\Cal T}_n(\eta)}{ \frac{|f^{\prime}(0)|}{\|f\|_{L_2{[0,1]}}}} \geq 
(1 + \varepsilon_n) \, 3^{-1/2} \, n^3$$ 
for every $n \in {\Bbb N}$ and for every $\eta > 0$, where $\varepsilon_n$ 
(determined exactly in the proof) is a quantity tending to $0$ an $n$ tends to $\infty$.  
\endproclaim

\medskip

\head 6. New Results: Markov-Nikolskii-type inequalities for $f \in {\Cal E}_n$ with nonnegative 
exponents \endhead

Our next theorem records how large $|f^{\prime}(0)|$ can be if $\|f\|_{L_2[0,1]} = 1$ for 
exponential sums $f \in {\Cal E}_n$ with nonnegative exponents.

\proclaim{Theorem 6.1}
We have
$$|f^{\prime}(0)| \leq (1 + \varepsilon_n) \, 3^{-1/2} \, n^3 \|f\|_{L_2[0,1]}$$
for every $f \in {\Cal E}_n^+$ of the form
$$f(t) = \sum_{j=1}^n{a_j e^{\lambda_j t}}\,, \qquad a_j \in {\Bbb R}\,, \enskip 0 \leq \lambda_1 < \lambda_2 < \cdots < \lambda_n\,,$$ 
where the quantity $\varepsilon_n$ (determined exactly in the proof) tends to $0$ an $n$ tends to $\infty$. 
This inequality is sharp.
\endproclaim  

\medskip

\head 7. New Results: A pointwise Nikolskii-type inequality for ${\Cal E}_n$ \endhead

The upper bound of the theorem below follows from Lemma 10.5 proved in [2]. We couple 
this upper bound with a a matching lower bound.  

\proclaim{Theorem 7.1} We have
$$\left( \frac{(n-2)\log 2}{32 \, \min \{y-a,b-y\}} \right)^{1/2} \leq 
\sup_{f \in {\Cal E}_n}{\frac{|f(y)|}{\|f\|_{L_2[a,b]}}} \leq 
\left( \frac{2n}{\min\{{y-a,b-y\}}} \right)^{1/2}$$
for every $y \in (a,b)$ 
\endproclaim

The theorem below shows a lower bound for 
$$\sup_{f \in {\Cal T}_n}{\frac{|f(y)|}{\|f\|_{L_2[a,b]}}}\,.$$
However, there is a gap between the lower bound of Theorem 7.2 and the upper bound of Theorem 7.1.

\proclaim{Theorem 7.2} There is an absolute constant $c > 0$ such that
$$c \, \min \left \{ \frac{n^{1/2}}{\left( \min \{y-a,b-y\} \right)^{1/4}}\,, \enskip \frac{n}{(b-a)^{1/2}} 
\right \} \leq \sup_{f \in {\Cal T}_n(\varepsilon)}{\frac{|f(y)|}{\|f\|_{L_2[a,b]}}}$$
for every $\varepsilon > 0$ and for every $y \in [a,b]$. 
\endproclaim

\medskip

\head 8. New Results: An essentially sharp Bernstein-type inequality for ${\Cal T}_n$ \endhead

Our next theorem may be viewed as an essentially sharp (up to the constant $2e$) Bernstein type inequality
for all $f \in {\Cal T}_n$ at least in the middle of the interval $[-1,1]$.

\proclaim{Theorem 8.1} We have
$$|f^{\prime}(0)| \leq (\lambda + 2e(n + 1)) \, \|f\|_{[-1,1]}$$
for every $f \in {\Cal T}_n$ of the form (1.6), where 
$$\lambda := \max_{1 \leq j \leq n} {|\lambda_j|}\,, \qquad \lambda_1 < \lambda_2 < \cdots < \lambda_n\,. \tag 8.1$$
This inequality is sharp up to the factor $2e$. Namely, for every real number $\lambda > 0$ and integer $n \geq 1$ 
there is an $f \in {\Cal T}_n$ of the form (1.6) with (2.1) such that 
$$|f^{\prime}(0)| \geq \frac 14 \, (\lambda + n - 3) \, \|f\|_{[-1,1]}\,.$$ 
\endproclaim

\medskip

\head 9. New Results: Markov-type inequality for ${\Cal T}_n$ \endhead 
The next theorem offers a Markov-type inequality for all $f \in {\Cal T}_n$ on $[0,1]$. 

\proclaim{Theorem 9.1} Let $\varepsilon_n$ be the same as in Theorems 3.1 and 2.14. We have 
$$|f^{\prime}(0)| \leq 
(1 + \varepsilon_n) \, \left( 27n^5 + \sum_{k=1}^n{\lambda_k^2} \right)^{1/2}\|f\|_{[0,1]}\,,$$
and
$$\|f^{\prime}\|_{[0,1]} \leq 
(1 + \varepsilon_n) \, \left( 108n^5 + \sum_{k=1}^n{\lambda_k^2} \right)^{1/2}\|f\|_{[0,1]}\,,$$
for every $f \in {\Cal T}_n$ of the form (1.6) 
\endproclaim
 
\proclaim{Theorem 9.2} We have 
$$\sup_{f \in {\Cal T}_n(\eta)}{ \frac{|f^{\prime}(0)|}{\|f\|_{[0,1]}}} \geq 2(n-1)^2$$
for every $n \in {\Bbb N}$ and for every $\eta > 0$. 
\endproclaim

\medskip

\head 10. New Results: Markov-type inequalities for ${\Cal E}_n$ and ${\Cal T}_n$ in $L_2[0,\infty)$ with the 
Laguerre weight \endhead

In this section we use the norm
$$\|f\| := \left( \int_0^\infty{|f(t)|^2e^{-t}\,dt} \right)^{1/2}\,.$$
Our first result extends Lubinsky's Theorem 3.2 in [22] to the case when 
the exponents are not necessarily purely imaginary.

\proclaim{Theorem 10.1}
We have
$$\|f^{\prime}\| \leq \left( \max_{1 \leq j \leq n}{|\lambda_j|} +
\left( \sum_{j=1}^n{(1 - 2{\text {\rm Re}}(\lambda_j)) \sum_{k=j+1}^n{(1 - 2{\text {\rm Re}}(\lambda_k))}}\right)^{1/2} \right) \|f\|$$
for every $f \in {\Cal E}_n$ of the form
$$f(t) = \sum_{j=1}^n{a_j e^{\lambda_jt}}, \quad a_j, \lambda_j \in {\Bbb C}\,, \enskip \text {\rm Re}(\lambda_j) < 1/2\,.$$
\endproclaim

The theorem below recaptures Lubinsky's Theorem 3.2 in [22].

\proclaim{Theorem 10.2}
We have
$$\|f^{\prime}\| \leq \left( \max_{1 \leq j \leq n}{|\lambda_j|} + \left( \frac{n(n-1)}{2} \right)^{1/2} \right) \|f\|$$
for every $f \in {\Cal T}_n$ of the form (1.6).
\endproclaim

\medskip

\head 11. New Results: Markov-type inequalities for ${\Cal E}_n^-$ in $L_2[0,\infty)$ \endhead

In this section we use the norm
$$\|f\| := \left( \int_0^\infty{|f(t)|^2 \,dt} \right)^{1/2}\,.$$

Our only result in this section is a version of Theorem 10.1, a Markov-type inequality for 
${\Cal E}_n^-$ in $L_2[0,\infty)$ without a weight.

\proclaim{Theorem 11.1}
We have
$$\|f^{\prime}\| \leq \left( \frac 12 + \max_{1 \leq j \leq n}{\left|\lambda_j + \frac 12 \right|} +
2 \left( \sum_{j=1}^n{{\text {\rm Re}}(\lambda_j) \sum_{k=j+1}^n{{\text {\rm Re}}(\lambda_k)}}\right)^{1/2} \right) \|f\|$$
for every $f \in {\Cal E}_n^-$ of the form
$$f(t) = \sum_{j=1}^n{a_j e^{\lambda_jt}}, \quad a_j, \lambda_j \in {\Bbb C}\,, \enskip \text {\rm Re}(\lambda_j) < 0\,.$$
\endproclaim

\medskip 

\head 12. Lemmas \endhead

Our first lemma is due to Tur\'an. See E.6 b] on page 297 of [3]. In fact, this inequality plays a central role 
in Tur\'an's book [29] as well.

\proclaim{Lemma 12.1} 
We have
$$|g(0)| \leq \left( \frac{2e(\alpha+\beta)}{\beta} \right)^n \|g\|_{[\alpha,\alpha+\beta]}\,, 
\qquad g \in {\Cal E}_n^+\,,$$  
for every $\alpha > 0$ and $\beta > 0$.
\endproclaim

In fact, we will need the following consequence of Lemma 12.2.

\proclaim{Lemma 12.2}
We have
$$|f(t)| \leq \left( \frac{2e(t-a)}{d} \right)^n \|f\|_{[a,a+d]} \leq 
\left( \frac{2et}{d} \right)^n \|f\|_{[a,a+d]}\,, \qquad f \in {\Cal E}_n^-\,,$$
for every $a > 0$, $d > 0$ and $t \geq a+d$.
\endproclaim

\demo{Proof of Lemma 12.2}
Let $f \in {\Cal E}_n^-$. Let $g \in {\Cal E}_n^+$ be defined by $g(x) := f(t-x)$.
Associated with $a > 0$, $d > 0$, $t \geq a+d$ we define $\alpha := t-(a+d)$, $\beta := d$. 
Applying Lemma 12.1 with $g \in {\Cal E}_n^+$ we get
$$\split |f(t)| = |g(0)| \leq & \, 
\left( \frac{2e(\alpha+\beta)}{\beta} \right)^n \|g\|_{[\alpha,\alpha+\beta]} \cr 
& \, = \left( \frac{2e(t-a)}{\beta} \right)^n \|f\|_{[a,a+d]} 
\leq \left( \frac{2et}{d} \right)^n \|f\|_{[a,a+d]} \cr \endsplit$$
\qed \enddemo

Our next lemma states the first inequality of part c] of E.2 coupled with 
part d] of E.2 on page 286 of [3]. See also Corollary 3.3 in [8]. 

\proclaim{Lemma 12.3}
We have 
$$\frac{|y^{1/2}P(y)|}{\|P\|_{L_2[0,1]}} \leq \left( \sum_{j=1}^n{(1 + 2{\text {\rm Re}}(\lambda_j))} \right)^{1/2}$$
for every M\"untz polynomial $0 \neq P$ of the form 
$$P(x) = \sum_{j=1}^n{a_j x^{\lambda_j}}, \quad a_j, \lambda_j \in {\Bbb C}\,, \enskip \text {\rm Re}(\lambda_j) > -1/2\,,$$
and for every $y \in [0,1]$. This inequality is sharp when $y=1$.
\endproclaim

Using the substitution $x = e^{-t}$ Lemma 12.3 implies the following.

\proclaim{Lemma 12.4}
We have 
$$|f(0)| \leq \left( \sum_{j=1}^n{(1 - 2{\text {\rm Re}}(\lambda_j))} \right)^{1/2} \, \|f(t) e^{-t/2}\|_{L_2[0,\infty)}$$
for every $f \in {\Cal E}_n$ of the form 
$$f(t) = \sum_{j=1}^n{a_je^{\lambda_jt}}\, \quad a_j \in {\Bbb C}\,, \enskip \text{\rm Re}(\lambda_j) < 1/2\,.$$ 
This inequality is sharp. 
\endproclaim

The next lemma is from [7].

\proclaim{Lemma 12.5}
We have
$$|f(y)| \leq \left( \frac{n}{\delta} \right)^{1/2} \|f\|_{L_2[y-\delta,y+\delta]}\,, \qquad f \in {\Cal E}_n\,,$$
for every $y \in {\Bbb R}$ and $\delta > 0$.
\endproclaim

Our next lemma states the second inequality of part c] of E.2 coupled with part d] of E.2 
on page 286 of [3]. See also Corollary 3.3 in [8].

\proclaim{Lemma 12.6}
We have
$$\frac{|y^{3/2}P^{\prime}(y)|}{\|P\|_{L_2[0,1]}} \leq \left( \sum_{k=1}^n{(1 + 2{\text {\rm Re}}(\lambda_k)) \, 
\Big| \lambda_k + \sum_{j=1}^{k-1}{(1 + 2{\text {\rm Re}}(\lambda_j))} \Big|^2} \right)^{1/2}$$
for every $y \in [0,1]$ and for every M\"untz polynomial $0 \neq P$ of the form
$$P(x) = \sum_{j=1}^n{a_j x^{\lambda_j}}, \quad a_j, \lambda_j \in {\Bbb C}\,, \enskip \text {\rm Re}(\lambda_j) > -1/2\,.$$
This inequality is sharp when $y=1$.
\endproclaim

Using the substitution $x = e^{-t}$ Lemma 12.6 implies the following.

\proclaim{Lemma 12.7}
We have
$$\frac{|f^{\prime}(0)|}{\|f(t)e^{-t/2}\|}_{L_2[0,\infty)} \leq \left( \sum_{k=1}^n{(1 + 2{\text {\rm Re}}(\lambda_k)) \, 
\Big| \lambda_k + \sum_{j=1}^{k-1}{(1 + 2{\text {\rm Re}}(\lambda_j))} \Big|^2} \right)^{1/2}$$
for every exponential sums $0 \neq f$ of the form
$$f(t) = \sum_{j=1}^n{a_j e^{i\lambda_j t}}, \quad a_j, \lambda_j \in {\Bbb C}\,, \enskip \text {\rm Re}(\lambda_j) < 1/2\,.$$
This inequality is sharp.
\endproclaim

The heart of the proof of our Theorem 4.1 is the following pair of
comparison lemmas. The proof of the next couple of lemmas is
based on basic properties of Descartes systems, in particular on
Descartes' Rule of Sign, and on a technique used earlier by
P.W. Smith and Pinkus. Lorentz ascribes this result to Pinkus,
although it was P.W. Smith [26] who published it. I have
learned about the the method of proofs of these lemmas from
Peter Borwein, who also ascribes it to Pinkus. This is the proof we 
present in [16]. Section 3.2 of [1], for instance, gives an introduction 
to Descartes systems. Descartes' Rule of Signs is stated and proved on 
page 102 of [3].

\proclaim{Lemma 12.8} Let $\Delta_n := \{\delta_0 < \delta_1 < \cdots < \delta_n\}$ and
$\Gamma_n := \{\gamma_0 < \gamma_1 < \cdots < \gamma_n\}$
be sets of real numbers satisfying $\delta_j \leq \gamma_j$ for each $j=0,1,\ldots,n$.
Let $a,b,c \in {\Bbb R}$, $a < b \leq c$. Let $w$ be a not identically $0$, continuous function
defined on $[a,b]$. Let $q \in (0,\infty]$. Then
$$\sup \left\{\frac{|P(c)|}{\|Pw\|_{L_q[a,b]}}: \enskip 0 \neq P \in E(\Delta_n)\right\}
\leq \sup \left\{\frac{|P(c)|}{\|Pw\|_{L_q[a,b]}}: \enskip 0 \neq P \in E(\Gamma_n)\right\}\,.$$
Under the additional assumption $\delta_n \geq 0$ we also have
$$\sup \left\{\frac{|P^{\prime}(c)|}{\|Pw\|_{L_q[a,b]}}: \enskip 0 \neq P \in E(\Delta_n)\right\}
\leq \sup \left\{\frac{|P^{\prime}(c)|}{\|Pw\|_{L_q[a,b]}}: \enskip
0 \neq P \in E(\Gamma_n)\right\}\,.$$
\endproclaim

\proclaim{Lemma 12.9} Let $\Delta_n := \{\delta_0 < \delta_1 < \cdots < \delta_n\}$ and
$\Gamma_n := \{\gamma_0 < \gamma_1 < \cdots < \gamma_n\}$
be sets of real numbers satisfying $\delta_j \leq \gamma_j$ for each $j=0,1,\ldots,n$.
Let $a,b,c \in {\Bbb R}$, $c \leq a < b$. Let $w$ be a not identically $0$, continuous function
defined on $[a,b]$. Let  $q \in (0,\infty]$. Then
$$\sup \left\{\frac{|P(c)|}{\|Pw\|_{L_q[a,b]}}: \enskip  0 \neq P \in E(\Delta_n)\right\}
\geq \sup \left\{\frac{|P(c)|}{\|Pw\|_{L_q[a,b]}}: \enskip 0 \neq P \in E(\Gamma_n)\right\}\,.$$
Under the additional assumption $\gamma_0 \leq 0$ we also have
$$\sup \left\{\frac{|Q^{\prime}(c)|}{\|Qw\|_{L_q[a,b]}}: \enskip  0 \neq Q \in E(\Delta_n)\right\}
\geq \sup \left\{\frac{|Q^{\prime}(c)|}{\|Qw\|_{L_q[a,b]}}: \enskip 0 \neq Q \in E(\Gamma_n)\right\}\,.$$
\endproclaim

An entire function $f$ is said to be of exponential type $\tau$ if for any $\varepsilon > 0$
there exists a constant $k(\varepsilon)$  such that 
$|f (z)| \leq k(\varepsilon)e^{(\tau + \varepsilon)|z|}$
for all $z \in {\Bbb C}$.
The following inequality may be found on p. 102 of [1] and is known as Bernstein's inequality.
See also [2] and [13]. It can be viewed as an extension of Bernstein's (trigonometric) polynomial 
inequality (see p. 232 of [3], for instance) to entire functions of exponential type bounded on 
the real axis.

\proclaim{Lemma 12.10 (Bernstein's inequality)}
Let $f$ be an entire function of exponential type $\tau > 0$ bounded on ${\Bbb R}$.
Then
$$\sup_{t \in {\Bbb R}}{|f^\prime(t)|} \leq \tau \sup_{t \in {\Bbb R}}{|f(t)|}\,.$$
\endproclaim

The reader may find another proof of the above Bernstein's inequality in [25, pp. 512--514],
where it is also shown that an entire function $f$ of exponential type $\tau$ satisfying
$$|f^\prime(t_0)| = \tau \sup_{t \in {\Bbb R}} |f(t)|$$
at some point $t_0 \in {\Bbb R}$ is of the form
$$f(z) = ae^{i\tau z} + be^{-i\tau z}\,, \qquad a \in {\Bbb C}, \quad b \in {\Bbb C}, \qquad
|a| + |b| = \sup_{t \in {\Bbb R}}{|f(t)|}\,.$$

Our next lemma is stated as Theorem 6.1.5 on page 282 of [3]. See also Theorem 3.4 in [8].

\proclaim{Lemma 12.11}
We have 
$$\|xP^{\prime}(x)\|_{L_2[0,1]} \leq \left( \sum_{j=1}^n{|\lambda_j|^2} + \sum_{j=1}^n{\sum_{j=1}^n{(1 + 2{\text {\rm Re}}(\lambda_j))
\sum_{k=j+1}^n}{(1 + 2{\text {\rm Re}}(\lambda_k))}}\right)^{1/2}\|P\|_{L_2[0,1]}$$
for every M\"untz polynomial $0 \neq P$ of the form
$$P(x) = \sum_{j=1}^n{a_j x^{\lambda_j}}, \quad a_j, \lambda_j \in {\Bbb C}\,, \enskip \text {\rm Re}(\lambda_j) > -1/2\,. \tag 12.1$$
\endproclaim

In fact, a simple change in the proof (in either references) gives the following.

\proclaim{Lemma 12.12}
We have
$$\|xP^{\prime}(x)\|_{L_2[0,1]}
\leq \left( \max_{1 \leq j \leq n}{|\lambda_j|} + 
\left( \sum_{j=1}^n{(1 + 2{\text {\rm Re}}(\lambda_j)) \sum_{k=j+1}^n{(1 + 2{\text {\rm Re}}(\lambda_k))}} \right)^{1/2} \right) 
\|P\|_{L_2[0,1]}$$
for every M\"untz polynomial $P$ of the form 12.1.
\endproclaim

\demo{Proof of Lemma 12.12}
Let $P$ be a M\"untz polynomial of the form (12.1). Then 
$$P(x) = \sum_{k=1}^n{a_kL_k^*}\,, \qquad a_k \in {\Bbb C}\,,$$
where 
$$L_k^* \in \text{span}\{x^{\lambda_1}, x^{\lambda_2}, \ldots, x^{\lambda_k}\}$$ 
denotes the $k$th orthonormal M\"untz-Legendre polynomials on $[0,1]$ associated with 
$$\text{span}\{x^{\lambda_1}, x^{\lambda_2}, \ldots, x^{\lambda_n}\}\,,$$
introduced in Section 3.4 of [3] (the spans here are taken over ${\Bbb C}$. 
Without loss of generality we may assume that
$$\|P\|_{L_2[0,1]} = \sum_{k=1}^n{|a_k|^2} = 1\,. \tag 12.2$$
As it is observed on page 283 of [3], we have
$$xP^\prime(x) = \sum_{j=1}^n{\left( a_j\lambda_j + 
\sqrt{1 + 2{\text {\rm Re}}(\lambda_j)} \sum_{k=j+1}^n{a_k \sqrt{1 + 2{\text {\rm Re}}(\lambda_k))}} \right) L_j^*(x)} \,$$
Hence
$$\|xP^{\prime}(x)\|_{L_2[0,1]} \leq \|R\|_{L_2[0,1]} + \|S\|_{L_2[0,1]}\,, \tag 12.3$$
where
$$R(x) := \sum_{j=1}^n{a_j \lambda_j L_j^*}$$
and 
$$S(x) := \sum_{j=1}^n{\left( \sqrt{1 + 2{\text {\rm Re}}(\lambda_j)} 
\sum_{k=j+1}^n{a_k \sqrt{1 + 2{\text {\rm Re}}(\lambda_k))}}\right) L_j^*(x)}\,.$$
Using the orthonormality of $\{L_j^*, j=1,2,\ldots,n\}$ on $[0,1]$ and then recalling (12.2), we can deduce that
$$\split \|R\|_{L_2[0,1]} = & \, \left( \sum_{j=1}^n{|a_j \lambda_j|^2} \right)^{1/2} \leq  
\max_{1 \leq j \leq n}{|\lambda_j|} \left( \sum_{j=1}^n{|a_j|^2} \right)^{1/2} \cr 
\leq & \, \max_{1 \leq j \leq n}{|\lambda_j|}\,. \cr \endsplit \tag 12.4$$
Further, combining the orthonormality of $\{L_j^*, j=1,2,\ldots,n\}$ on $[0,1]$ with applications of the 
Cauchy-Schwarz inequality to each term of the first sum and then recalling (12.2) we obtain that
$$\split \|S\|_{L_2[0,1]}^2 
= & \, \sum_{j=1}^n{(1 + 2{\text {\rm Re}}(\lambda_j)) 
\left| \sum_{k=j+1}^n{a_k \sqrt{1 + 2{\text {\rm Re}}(\lambda_k))}} \right|^2} \cr
\leq & \,\sum_{j=1}^n{(1 + 2{\text {\rm Re}}(\lambda_j)) \sum_{k=j+1}^n{(1 + 2{\text {\rm Re}}(\lambda_k))}} \cr
\endsplit \tag 12.5$$
The lemma now follows from (12.3), (12.4) and (12.5).
\qed \enddemo

\medskip

\head 13. Proofs of the new results \endhead 

\demo{Proof of Theorem 2.1}
Let $f \in {\Cal T}_n$. By Lemma 12.5 we have
$$\|f\|_{[n,7n]} \leq  \|f\|_{L_2[0,8n]}\,.$$
Combining this with Lemma 12.2 we get
$$\split |f(t)|^2 e^{-t} \leq  & \, \left( \left( \frac{2et}{6n} \right)^{2n} \|f\|_{[n,7n]}^2 \right) e^{-t} 
\leq \left( \left( \frac{et}{3n} \right)^{2n}  \|f\|_{L_2[0,8n]}^2 \right) e^{-t} \cr 
\leq & \, e^{-t/2} \|f\|_{L_2[0,8n]}^2\,, \qquad t \geq 8n\,. \cr \endsplit$$
Here we used the fact that 
$$h(t) := \left( \frac{et}{3n} \right)^{2n} e^{-t/2}$$
is decreasing on the interval $[8n,\infty)$, hence
$$\split \left( \frac{et}{3n} \right)^{2n} e^{-t} \leq & \, \left( \left( \frac{et}{3n} \right)^{2n} e^{-t/2} \right) e^{-t/2}  \leq \left( \frac{8e}{3} \right)^{2n} e^{-4n} e^{-t/2} \cr 
\leq & \left( \frac{(8/3)^2 e^2}{e^4} \right)^n \leq \, e^{-t/2}\,, \quad t \geq 8n\,. \cr \endsplit$$
Hence
$$\int_{8n}^{\infty}{|f(t)|^2 e^{-t} \, dt} \leq \left( \int_{8n}^{\infty}{e^{-t/2} \, dt} \right) \, \|f\|_{L_2[0,8n]}^2 \leq 2e^{-2n} \int_0^{8n}{|f(t)|^2 e^{-t/4} \, dt}\,.$$
This implies that
$$\int_0^{\infty}{|f(t)|^2 e^{-t} \, dt} \leq (1+2e^{-2n}) \, \int_0^{8n}{|f(t)|^2 e^{-nt/4} \, dt}\,.$$
Combining this with Lemma 12.4 we get
$$|f(0)| \leq (1+2e^{-2n})^{1/2} \, n^{1/2} \, \|f(t)e^{-t/8}\|_{L_2[0,8n]}$$ 
Transforming this inequality linearly from the interval $[0,8n]$ to the interval $[0,1]$, we get the first statement of the theorem. 

The second statement of the theorem follows from the second statement of Lemma 12.4. 
Indeed, for every fixed $\lambda_1 < \lambda_2 < \cdots < \lambda_n$ there is a $0 \neq f \in {\Cal T}_n$ of the form (1.6) 
such that 
$$|f(0)| \geq  n^{1/2} \, \|f(t) e^{-t/2}\|_{L_2[0,\infty)} > n^{1/2} \, \|f(t) e^{-t/2}\|_{L_2[0,8n]}\,.$$ 
Transforming this inequality linearly from the interval $[0,8n]$ to the interval $[0,1]$, we get the second statement of the theorem.
\qed \enddemo

\demo{Proof of Theorem 2.2}
Let $\gamma_0 := 2 + \log 4 < \gamma \leq 4$ and $\delta := (\gamma - \gamma_0)/8 < 1/8$. 
Observe that $\gamma_0 < \gamma \leq 4$ implies that $0 < \delta < 1/8$ and hence
$$\gamma - 2\delta \geq \gamma_0 - 2\delta \gamma_0 - 1/4 > 2\,.$$
Combining this with the Mean Value Theorem we obtain
$$\log \gamma - \log (\gamma - 2\delta) < 2\delta \,\frac{1}{\gamma-2\delta} \geq 
2\delta \, \frac 12 = \delta\,.$$
Therefore 
$$\split 2 + \log 4 + 2\log \frac{\gamma}{\gamma-2\delta} - \gamma + \gamma \delta = & \, 
(\gamma_0 - \gamma) + 2(\log \gamma - \log (\gamma - 2\delta)) + \gamma \delta \cr 
< & \, -8\delta + 2\delta + 4 \delta = -2\delta < 0\,, \cr \endsplit$$
hence
$$4e^2 \left( \frac{\gamma}{\gamma-2\delta} \right)^2e^{\gamma (\delta - 1)} \leq 1\,. \tag 13.1$$   
Let $f \in {\Cal T}_n$. By Lemma 12.5 we have
$$\|f\|_{[\delta n, (\gamma - \delta)n]}^2 \leq  \delta^{-1} \|f\|_{L_2[0,\gamma n]}^2\,.$$
Combining this with Lemma 12.2 we get
$$\split |f(t)|^2 e^{-t} \leq  & \, \left( \left( \frac{2e(t-\delta n)}{(\gamma - 2\delta)n} \right)^{2n} \|f\|_{[\delta n,(\gamma - \delta)n]}^2 \right) e^{-t} \cr 
\leq & \, \delta^{-1}\left( \left( \frac{2et}{(\gamma - 2\delta)n} \right)^{2n}  \|f\|_{L_2[0,\gamma n]}^2 \right) e^{-t} \cr 
\leq & \, \delta^{-1} e^{-\delta t} \|f\|_{L_2[0,\gamma n]}^2\,, \qquad t \geq \gamma n\,. \cr \endsplit$$
Here we used the fact that
$$h(t) := \left( \frac{2et}{(\gamma - 2\delta)n} \right)^{2n} e^{(\delta -1)t}$$
is decreasing on the interval $[\gamma,\infty) \subset [2(1-\delta)^{-1},\infty)$, which, together with (13.1) yields
$$\split \left( \frac{2et}{(\gamma - 2\delta)n} \right)^{2n} e^{-t} \leq & \, 
\left( \left( \frac{2e\gamma}{\gamma - 2\delta)} \right)^{2n} e^{(\delta -1)t} \right) e^{-\delta t}
\left( \left( \frac{2e\gamma}{\gamma - 2\delta)} \right)^{2n} e^{\gamma(\delta -1)n} \right) e^{-\delta t} \cr 
\leq & \left( 4e^2 \left( \frac{\gamma}{\gamma-2\delta} \right)^2e^{\gamma (\delta - 1)} \right)^n \, e^{-\delta t}\,, 
\leq e^{-\delta t} \quad t \geq \gamma n\,, \cr \endsplit$$
Hence
$$\split \int_{\gamma n}^{\infty}{|f(t)|^2 e^{-t} \, dt} \leq & \, \delta^{-1} \left( \int_{\gamma n}^{\infty}{e^{-\delta t} \, dt} \right) \, \|f\|_{L_2[0,\gamma n]}^2 \cr  
\leq & \, \delta^{-1} \delta^{-1} e^{-\delta \gamma n} \int_0^{\gamma n}{|f(t)|^2 \, dt}\,. \cr \endsplit$$
This implies that
$$\int_0^{\infty}{|f(t)|^2 e^{-t} \, dt} \leq (1+ \delta^{-2} e^{-\delta \gamma n}) \, \int_0^{\gamma n}{|f(t)|^2\, dt}\,.$$
Combining this with Lemma 12.3 we get
$$|f(0)| \leq n^{1/2} \, \|f\|_{L_2[0,\infty)} \leq 
(1+ \delta^{-2} e^{-\delta \gamma n})^{1/2}\, n^{1/2} \, \|f\|_{L_2[0,\gamma n]}\,.$$
Transforming this inequality linearly from the interval $[0,\gamma n]$ to the interval $[0,1]$, we get the theorem.
\qed \enddemo

\demo{Proof of Theorem 2.3}
Let $y \in [-1,1]$. Transforming the inequality of Theorem 2.1 (with the constant $\pi/2$ rather than $(\gamma + \varepsilon_n)^{1/2}$) linearly to the intervals 
$[0,y]$ and $[y,1]$
$$y \, |f(y)|^2 \leq \left( \frac {\pi n}{2} \right)^2 \int_{[0,y]}{|f(t)|^2 \, dt}$$
and
$$(1-y) \, |f(y)|^2 \leq \left( \frac {\pi n}{2} \right)^2 \int_{[y,1]}{|f(t)|^2 \, dt}\,.$$    
Adding these, we conclude that
$$|f(y)|^2 \leq \left( \frac {\pi n}{2} \right)^2 \int_{[0,1]}{|f(t)|^2 \, dt}\,,$$  
and the theorem follows.
\qed \enddemo

\demo{Proof of Theorem 2.4}
Let $f \in {\Cal T}_n$ and $q \in (0,2]$. Using Theorem 2.3 we obtain
$$\split \|f\|_{[0,1]} \leq & \frac {\pi n}{2} \, \|f\|_{L_2[0,1]} =  
\frac {\pi n}{2} \, \left( \int_0^1{|f(t)|^2 \, dt} \right)^{1/2} \cr 
\leq & \, \frac {\pi n}{2} \, \, \left( \int_0^1{|f(t)|^q \|f\|_{[0,1]}^{2-q} \, dt} \right)^{1/2}\,, \cr
\endsplit$$
and hence
$$\|f\|_{[0,1]}^{q/2} \leq \frac {\pi n}{2} \, \|f\|_{L_q[0,1]}^{q/2}\,,$$
and the theorem follows.
\qed \enddemo

\demo{Proof of Theorem 2.5}
When $p = \infty$ and $q \in (0,2]$, the theorem follows from Theorem 2.4.
Let $0 < q < p < \infty$, $q \leq 2$, and $f \in {\Cal T}_n$.
Based on Theorem 2.4 the proof of the theorem is fairly routine. We have
$$\split \|f\|_{L_p[0,1]}^p  = & \, \int_{[0,1]}{|f(t)|^p \, dt} \leq \int_{[0,1]}{|f(t)|^q \|f\|_{[0,1]}^{p-q} \, dt} \cr 
\leq & \, \|f\|_{L_q[0,1]}^q \|f\|_{[0,1]}^{p-q} \leq \|f\|_{L_q[0,1]}^q \left( \frac {\pi n}{2} \right)^{(p-q)2/q} \|f\|_{L_q[0,1]}^{p-q} \cr
\leq & \, \left( \frac {\pi n}{2} \right)^{(p-q)2/q} \|f\|_{L_q[0,1]}^p\,, \cr \endsplit$$
and by taking the $p$th root of both sides the theorem follows.
\qed \enddemo

\demo{Proof of Theorem 2.6}
The remark following Theorem 7.17.1 on page 182 of [27] asserts that
$$\sup_{P \in {\Cal P}_n}{\frac{|P(1)|}{\|P\|_{L_2{[-1,1]}}}} = 
\sup_{P \in {\Cal P}_n}{\frac{\|P\|_{[-1,1]}}{\|P\|_{L_2{[-1,1]}}}} = 2^{-1/2}(n+1)\,,$$
where ${\Cal P}_n$ denotes the set of all algebraic polynomials of degree at most $n$ with real coefficients.
Combining this with the observation
$$t := \lim_{\varepsilon \rightarrow 0+}{\frac{e^{i\varepsilon t} - 1}{i\varepsilon}} \,, \tag 13.2$$
the theorem follows by a linear transformation from $[-1,1]$ to $[0,1]$.
\qed \enddemo

\demo{Proof of Theorem 2.7}
The guided exercise E.19 on page 413 of [3] shows that  
$$\sup_{P \in {\Cal P}_n}{\frac{|P(1)|}{\|P\|_{L_q{[-1,1]}}}} = \sup_{P \in {\Cal P}_n}{\frac{\|P\|_{[-1,1}]}{\|P\|_{L_q{[-1,1]}}}} 
\geq c^{1+1/q}(1+qn)^{2/q}$$
for every $q \in (0,\infty)$. Combining this with the observation (13.2) the theorem follows by a linear transformation from $[-1,1]$ to $[0,1]$.
\qed \enddemo

\demo{Proof of Theorem 2.9}
Let $q \in (2,\infty)$ and let $1/p := (q-2)/q$, that is, $1/p + 1/(q/2) = 1$.  
Using Theorem 2.1 and H\"older's inequality, we have
$$\split |f(0)|^2 \leq & \, (8 + \varepsilon_n) \, n^2 \int_0^1{|f(t)|^2e^{-nt} e^{-nt} \, dt} \cr  
\leq & \, (8 + \varepsilon_n) \, n^2 \left( \int_0^1{\big( |f(t)|^2 e^{-nt} \big)^{q/2} dt} \right)^{2/q} 
\left( \int_0^1{\big|e^{-nt}\big|^p\, dt} \right)^{1/p} \,, \cr \endsplit$$
hence
$$\split |f(0)| \leq & \, (8 + \varepsilon_n)^{1/2} \, n \|f(t)e^{-nt}\|_{L_q[0,1]} \left( \frac{1}{pn} \right)^{1/p} \cr 
\leq & \, (8 + \varepsilon_n)^{1/2} \, n \, \|f(t)e^{-nt}\|_{L_q[0,1]} \left( \frac{q-2}{2qn} \right)^{(q-2)/q} \cr
\leq & \, (8 + \varepsilon_n)^{1/2} \, c_q n^{1/2 + 1/q} \|f(t)e^{-nt}\|_{L_q[0,1]} \,. \cr \endsplit$$
\qed \enddemo

\demo{Proof of Theorem 2.10}
Let $y \in [0,1]$. Transforming the inequality of Theorem 2.9 linearly to the intervals $[0,y]$ and $[y,1]$, respectively, we obtain that
$$y \, |f(y)|^q \leq \left((8 + \varepsilon_n)^{1/2} c_q n^{1/2 + 1/q} \right)^q \, \int_0^y{|f(t)|^q \, dt}$$
and
$$(1-y) \, |f(y)|^q \leq  \left(( 8 + \varepsilon_n)^{1/2} c_q n^{1/2 + 1/q} \right)^q \, \int_y^1{|f(t)|^q \, dt}\,.$$
Adding these we conclude that
$$|f(y)|^q \leq \left( (8 + \varepsilon_n)^{1/2} c_q n^{1/2 + 1/q} \right)^q \, \int_0^1{|f(t)|^q \, dt}\,,$$
and the theorem follows.
\qed \enddemo

\demo{Proof of Theorem 3.1}
Let $f \in {\Cal E}_n^-$. Let $\delta := 1/91$ and $\eta := 1/90$. By Lemma 12.5 we have
$$\|f\|_{[\delta n, (2 -\delta)n]} \leq  \delta^{-1} \|f\|_{L_2[0,2n]}\,.$$
Combining this with Lemma 12.2 we get
$$\split |f(t)|^2 e^{-t} \leq  & \, \left( \left( \frac{2et}{(2-2\delta)n} \right)^{2n} 
\|f\|_{[\delta n,(2-\delta) n]}^2 \right) e^{-t} 
\leq \left( \left( \frac{2et}{(2-2\delta)n} \right)^{2n}  \delta^{-1} \|f\|_{L_2[0,2n]}^2 \right) e^{-t} \cr 
\leq & \,\delta^{-1} \left( \frac{2et}{(2-2\delta)n} \right)^{2n} e^{2n} e^{-t} \left( \int_0^{2n}{|f(x)|^2 e^{-x} \,dx} \right) \,, \qquad t \geq 2n\,. \cr \endsplit$$
Integrating on $[9n,\infty]$, we get
$$\split & \, \int_{9n}^{\infty}{|f(t)|^2 e^{-t} \, dt} 
\leq \, \delta^{-1} \left( \int_{9n}^{\infty}{\left( \frac{2et}{(2-2\delta)n} \right)^{2n} e^{2n} e^{-t} \, dt} \right) \left( \int_0^{2n}{|f(x)|^2 e^{-x} \, dx} \right) \cr
= & \, \delta^{-1} \left( \sup_{t \geq 9n}{\left( \frac{2et}{(2-2\delta)n} \right)^{2n} e^{2n} e^{(\eta-1)t} \, dt} \right) \left( \int_{9n}^{\infty}{e^{-\eta t} \, dt} \right)  
\left( \int_0^{2n}{|f(x)|^2 e^{-x} \,dx} \right) \cr
\leq & \, \delta^{-1} \left( \int_{9n}^{\infty}{e^{-\eta t} \, dt} \right) \left( \int_0^{2n}{|f(x)|^2 e^{-x} \, dx} \right) \cr 
\leq & \, \delta^{-1} \eta^{-1} e^{-9\eta n} \left( \int_{0}^{2n}{|f(x)|^2 e^{-x} \, dx} \right). \cr \endsplit \tag 13.3$$
Here we used the fact that
$$h(t) := \left( \frac{2et}{(2-2\delta)n} \right)^{2n} e^{2n} e^{(\eta -1)t}$$
is decreasing on the interval $[9n,\infty)$, hence recalling that $\delta := 1/91$ and $\eta = 1/90$, we have
$$\sup_{t \geq 9n}{h(t)}  \leq ((9.1)e)^{2n} e^{-(8.9)n} e^{2n} = e^{(2\log(9.1) + 2 - 8.9 + 2)n} \leq e^0 = 1\,.$$
It follows from (13.3) that
$$\int_{9n}^{\infty}{|f(t)|^2 e^{-t} \, dt} \leq \delta^{-1} \eta^{-1} e^{-9\eta n} \left( \int_0^{2n}{|f(x)|^2 e^{-x} \, dx} \right)\,,$$
hence
$$\int_0^{\infty}{|f(t)|^2 e^{-t} \, dt} \leq (1 + \delta^{-1} \eta^{-1} e^{-9\eta n}) \, \left( \int_0^{9n}{|f(x)|^2 e^{-x} \, dx} \right) \,.$$
\qed \enddemo

\demo{Proof of Theorem 3.2}
Let $f \in {\Cal T}_n$. Then, Lemma 12.4, yields that
$$|f(0)|^2 \leq n \, \int_0^{\infty} {|f(t)|^2 e^{-t} \, dt}\,.$$
Combining this with Theorem 3.1 we have
$$|f(0)|^2 \leq (1 + \varepsilon_n)^2 n \, \int_0^{9n} {|f(t)|^2 e^{-t} \, dt}\,.$$
Transforming this inequality from the interval $[0,9n]$ to the interval $[0,1]$, we obtain
$$|f(0)|^2 \leq  (1 + \varepsilon_n)^2 9n^2 \int_0^1{|f(u)|^2 \, e^{-9nu} \, du}\,.$$

The second statement of the theorem follows from the second statement of Lemma 12.4.
Indeed, for every fixed $\lambda_1 < \lambda_2 < \cdots < \lambda_n$ there is a $0 \neq f \in {\Cal T}_n$ of the form (1.6)
such that
$$|f(0)| =  n^{1/2} \, \|f(t) e^{-t/2}\|_{L_2[0,\infty)} > n^{1/2} \, \|f(t) e^{-t/2}\|_{L_2[0,9n]}\,.$$
Transforming this inequality linearly from the interval $[0,9n]$ to the interval $[0,1]$, we get the second statement of the theorem.
\qed \enddemo

\demo{Proof of Theorem 4.1}
Observe that
$$t = \lim_{\varepsilon \rightarrow 0+}{\frac{e^{\varepsilon t} - 1}{\varepsilon}}\,,$$
Hence it follows from Lemma 12.8 in a routine fashion that it is sufficient to prove 
the inequality only for polynomials $P \in {\Cal P}_{n-1}$, where ${\Cal P}_{n-1}$ denotes the set of 
all polynomials of degree at most $n-1$ with real coefficients, and this has been done 
in the proof of Theorem 2.6. The sharpness of the theorem also follows from the proof of 
Theorem 2.6.
\qed \enddemo

\demo{Proof of Theorem 5.1}
Let $f \in {\Cal T}_n$ be of the form (1.6), and let $g \in {\Cal T}_n$ be defined by 
$g(9nt) := f(t)$. By Theorem 3.1 we have 
$$\int_0^{\infty}{|g(t)|^2 e^{-t} \, dt} \leq (1 + \varepsilon_n)^2 \int_0^{9n}{|g(t)|^2 e^{-t} \, dt}\,.$$
Combining this with Lemma 12.7 we get
$$\split |f^{\prime}(0)| = & \, 9n|g^{\prime}(0)| \cr \leq & \, 9n 
\left( \sum_{k=1}^n{\left( \left( \frac{\lambda_k}{9n} \right)^2 + (k-1)^2 \right)} \right)^{1/2}
\|g(t)e^{-t/2}\|_{L_2[0,\infty)} \cr
\leq & \, 9n(1 + \varepsilon_n) \, 
\left( \sum_{k=1}^n{\left( \left( \frac{\lambda_k}{9n} \right)^2 + (k-1)^2 \right)} \right)^{1/2}
\|g(t)e^{-t/2}\|_{L_2[0,9n]} \cr
= & \, 9n(1 + \varepsilon_n) \, 
\left( \sum_{k=1}^n{\left( \left( \frac{\lambda_k}{9n} \right)^2 + (k-1)^2 \right)} \right)^{1/2}
3n^{1/2}\|f(u)e^{-9nu/2}\|_{L_2[0,1]} \,.\cr
\endsplit$$

The second statement of the theorem follows from the second statement of Lemma 12.7.
Indeed, for every fixed $\lambda_1 < \lambda_2 < \cdots < \lambda_n$ there is a $g \in {\Cal T}_n$ such that $f \in {\Cal T}_n$ defined by 
$g(9nt) := f(t)$ is of the form (1.6) and
$$\split |f^{\prime}(0)| = & \, 9n|g^{\prime}(0)| \cr  
= & \, 9n \, \left( \sum_{k=1}^n{\left( \left( \frac{\lambda_k}{9n} \right)^2 + (k-1)^2 \right)} \right)^{1/2} \, \|g(t) e^{-t/2}\|_{L_2[0,\infty)} \cr 
> & \, 9n \, \left( \sum_{k=1}^n{\left( \left( \frac{\lambda_k}{9n} \right)^2 + (k-1)^2 \right)} \right)^{1/2} \, 3n^{1/2}\|f(u)e^{-9nu/2}\|_{L_2[0,1]}\,. \cr \endsplit$$    
\qed \enddemo

\demo{Proof of Theorem 5.2}
Let $y \in [0,1]$. Transforming the inequality of Theorem 5.1 linearly to the intervals $[0,y]$ and $[y,1]$, respectively, we obtain that
$$y^3 \, |f^{\prime}(y)|^2 \leq \,  27^2(1 + \varepsilon_n)^2 \, n^3 
\left( \sum_{k=1}^n{\left( \left( \frac{y\lambda_k}{9n} \right)^2 + (k-1)^2 \right)} \right)^{1/2} \, 
\int_0^y{|f(t)|^2 \, dt}$$
and
$$(1-y)^3 \, |f^{\prime}(y)|^2 \leq \, 27^2(1 +\varepsilon_n)^2 \, n^3 
\left( \sum_{k=1}^n{\left( \left( \frac{(1-y)\lambda_k}{9n} \right)^2 + (k-1)^2 \right)} \right)^{1/2} 
\int_y^1{|f(t)|^2 \, dt}\,.$$
Using the first inequality above if $y \in [1/2,1]$ and the second inequality above if $y \in [1/2,1]$ we conclude that
$$|f^{\prime}(y)|^2 \leq 27^2 \, (1 +\varepsilon_n)^2 \, n^3 
\left( \sum_{k=1}^n{\left( 2\left( \frac{\lambda_k}{9n} \right)^2 + 8(k-1)^2 \right)} \right) \, 
\int_0^1{|f(t)|^2 \, dt}\,,$$
and the theorem follows.
\qed \enddemo

\demo{Proof of Theorem 5.3} 
Observe that
$$t = \lim_{\varepsilon \rightarrow 0+}{\frac{e^{i\varepsilon t} - 1}{i\varepsilon}}\,, \tag 13.4$$ 
hence 
$$\sup_{f \in {\Cal T}_n^c}{\frac{|f^{\prime}(0)|}{\|f\|_{L_2[0,1]}}} =  
\sup_{P \in {\Cal P}_{n-1}^c}{\frac{|P^{\prime}(0)|}{\|P\|_{L_2[0,1]}}} \tag 13.5$$
with an absolute constant $c>0$, where ${\Cal P}_n^c$ denotes the set of all polynomials of degree 
at most $n$ with complex coefficients. 
Let $P_n \in {\Cal P}_n$ be the $n$-th orthonormal Legendre polynomial on the interval $[0,1]$, that is,
$$\int_{0}^1{P_n(x)P_m(x) \, dx} = \, \delta_{n,m}\,.$$
where $\delta_{n,m} = 1$ if $n=m$ and $\delta_{n,m} = 0$ if $n \neq m$.
Recall that
$$P_k^{\prime}(0) = (-1)^k k(k+1)(2k+1)^{1/2}\,, \qquad k=0,1,\ldots\,. \tag 13.6$$
This can be seen by combining (4.21.7), (4.3.3), and (4.1.4) in [27] and by using a linear transformation 
from the interval $[-1,1]$ to the interval $[0,1]$.
As a consequence of orthonormality, the Cauchy-Schwarz inequality, and (13.6) it is well known 
(see E.2 on page 285 of [3], for instance) that
$$\split \sup_{P \in {\Cal P}_{n-1}}{\frac{|P^\prime(0)|}{\|P\|_{L_2[0,1]}}} = & \, 
\left( \sum_{k=0}^{n-1} {P_k^{\prime}(1)^2} \right)^{1/2} \cr
= & \, \left( \sum_{k=0}^{n-1} {k^2(k+1)^2(2k+1)} \right)^{1/2} = (1 + \varepsilon_n) \, 3^{-1/2} \, n^3 \,. \cr \endsplit \tag 13.7$$
Combining (13.5) and (13.7) gives the theorem.
\qed \enddemo

\demo{Proof of Theorem 6.1}
It follows from (13.4) and Lemma 12.9 in a routine fashion that it is sufficient to prove 
the inequality only for polynomials $P \in {\Cal P}_{n-1}$, where ${\Cal P}_{n-1}$ 
denotes the set of all polynomials of degree at most $n-1$ with real coefficients. 
Hence, combining (13.5) and (13.7) gives the theorem.
\qed \enddemo

\demo{Proof of Theorem 7.1}
The upper bound follows Lemma 12.5, see [7] for a proof. To see the lower bound we proceed as 
follows. Let $P_n \in {\Cal P}_n$ be the $n$-th orthonormal Legendre polynomial on the interval 
$[-1,1]$, that is, 
$$\int_{-1}^1{P_n(x)P_m(x) \, dx} = \, \delta_{n,m}\,,$$
where $\delta_{n,m} = 1$ if $n=m$ and $\delta_{n,m} = 0$ if $n \neq m$. Let
$$Q(x) = \sum_{k=0}^n{P_k(0)P_k(x)}\, \tag 13.8$$
Then
$$\|Q\|_{L_2[-1,1]}^2 = \sum_{k=0}^n{P_k(0)^2} \qquad \text{and} \qquad 
|Q(0)| = \sum_{k=0}^n{P_k(0)^2}\,, \tag 13.9$$
hence
$$\frac{|Q(0)|^2}{\|Q\|_{L_2[-1,1]}^2} = \sum_{k=0}^n{L_k(0)^2}\,.$$
It is well known (see p. 165 of [27], for example) that $P_k(0) = 0$ if $k$ is even, and
$$\split |P_k(0)|^2 = & \, \frac{2k+1}{2} 
\left( \frac 12 \right)^2 \left( \frac 34 \right)^2 \left( \frac 56 \right)^2  \cdots 
\left( \frac{k-3}{k-2} \right)^2 \left( \frac {k-1}{k} \right)^2 \cr 
\geq & \,  \left( \frac 12 \right)^2 \left( \frac 23 \, \frac 34 \right) 
\left( \frac 45 \, \frac 56 \right) \cdots 
\left( \frac {k-4}{k-3} \, \frac {k-3}{k-2} \right) \left( \frac {k-2}{k-1} \, 
\frac {k-1}{k} \right) \cr 
\geq & \, \frac {2k+1}{4k} \geq \frac 12 \cr \endsplit$$
if $k$ is odd.
Combining this with (13.8) and (13.9) gives
$$\frac{|Q(0)|^2}{\|Q\|_{L_2[-1,1]}^2} \geq  \frac{n-2}{4} \,.$$
Let $f(t) = Q(2e^{-t}-1)e^{-t/2}$. Then
$$\frac{|f(\log 2)|}{\|f\|_{L_2[0,\infty)}} = 
\frac{|Q(0)|}{2^{1/2} \, \|Q\|_{L_2[-1,1]}} \geq \frac{(n-2)^{1/2}}{8^{1/2}} \,.$$
Transforming the above inequality linearly from the interval $[0,\infty)$ to $[a,\infty)$ 
and $(-\infty,b]$, we get the the lower bound of the theorem. 
\qed \enddemo

\demo{Proof of Theorem 7.2}
Theorem 2.1 of [18] implies that there is an absolute constant $c > 0$ such that
$$c \, \min \left \{ \frac{n^{1/2}}{(1-y^2)^{1/4}}\,,\enskip  n \right \} \leq 
\sup_{P \in {\Cal P}_{n-1}} {\frac{|P(y)|}{\|P\|_{[-1,1]}}}\,,$$
for every $y \in [-1,1]$, where ${\Cal P}_{n-1}$ denotes the set of all algebraic polynomials of 
degree at most $n-1$ with real coefficients. Hence the theorem follows from (13.4).
\qed \enddemo

\demo{Proof of Theorem 8.1}
Let $f \in {\Cal T}_n$ be of the form (1.6) with (8.1). 
Let $m$ be an integer such that $n \leq 2m$. We define the entire function of type $\lambda = 2m$ by  
$$g(z) := f(z) \left( \frac{\sin z}{z} \right)^{2m}\,.$$ 
By Bernstein's inequality we have
$$|f^{\prime}(0)| = |g^{\prime}(0)| \leq (\lambda + 2m) \sup_{t \in {\Bbb R}}{|g(t)|}\,. \tag 13.10$$
Lemma 12.2 implies that
$$|g(t)| \leq \left( \frac{2et}{2e} \right)^n \|f\|_{[0,2e]} \left( \frac{|\sin t|}{t} \right)^{2m} 
\leq t^{n-2m} \|f\|_{[0,2e]} \leq \|f\|_{[0,2e]} \,, \qquad t \geq 2e\,, \tag 13.11$$
and as $|\sin t| \leq t$ for all $t \geq 0$ obviously 
$$|g(t)| \leq |f(t)|\,, \qquad t \in [0,2e]\,. \tag 13.12$$
Combining (13.11) and (13.12) we have 
$$\sup_{t \in [0,\infty)}{|g(t)|} \leq \|f\|_{[0,2e]}\,, \tag 13.13$$
and similarly 
$$\sup_{t \in (-\infty,0]}{|g(t)|} \leq \|f\|_{[-2e,0]}\,. \tag 13.14$$
Using (13.10), (13.13), and (13.14) we conclude
$$|f^{\prime}(0)| \leq (\lambda + 2m) \, \|f\|_{[-2e,2e]}\,.$$
Transforming the above inequality linearly from the interval $[-2e,2e]$ to the interval $[-1,1]$, 
and choosing $m$ so that $n=2m$ in $n$ is even, and $n+1=2m$ if $n$ is odd, we get the upper bound of the theorem. 
To see the sharpness of the upper bound up to the factor $2e$, we pick $f(t) := \sin {\lambda t}$ if $\lambda \geq n \geq 2$, and  
$f(t) = T_m(\varepsilon^{-1}\sin(\varepsilon t))$ with a sufficiently small $\varepsilon > 0$, 
where $T_m$ is the Chebyshev polynomial of degree $m$ defined by $T_m(\cos \theta) = \cos(m\theta)$, 
$\theta \in [0,2\pi)$, and $m$ is the largest odd integer such that $2m+1 \leq n$.     
\qed \enddemo

\demo{Proof of Theorem 9.1}
Let $y \in [0,1]$. Let $f \in {\Cal T}_n$ be of the form (1.6). Transforming the inequality of 
Theorem 5.1 linearly from the interval $[0,1]$ 
to the intervals $[0,y]$ and $[y,1]$, respectively, we obtain that
$$\split y^3|f^{\prime}(y)|^2 \leq & \,  27^2 \, (1 + \varepsilon_n)^2 \, n^3 
\left( \sum_{k=1}^n{\left( \left( \frac{\lambda_k}{9n} \right)^2 + (k-1)^2\right)} \right) 
\int_0^y{|f(u)|^2 e^{-9n(y-u)/y} \, du} \cr
\leq & \, 27^2 \, (1+\varepsilon_n)^2 \, n^3 
\left( \sum_{k=1}^n{\left( \left( \frac{y\lambda_k}{9n} \right)^2 + (k-1)^2 \right)} \right) \, 
\frac{y}{9n} \, \|f\|_{[0,y]}^2\,, \cr \endsplit$$
and
$$\split (1-y)^3 & |f^{\prime}(y)|^2 \cr  \leq & \, (1 + \varepsilon_n)^2 \, n^3
\left( \sum_{k=1}^n{\left( \left( \frac{(1-y)\lambda_k}{n} \right)^2 + (k-1)^2\right)} \right)
\int_0^y{|f(u)|^2 e^{-9n(y-u)/(1-y)} \, du} \cr 
\leq & \, 27^2 \, (1 + \varepsilon_n)^2 \, n^3 
\left( \sum_{k=1}^n{\left( \left( \frac{(1-y)\lambda_k}{9n} \right)^2 + (k-1)^2\right)} \right) \, 
\frac{1-y}{9n} \|f\|_{[y,1]}^2\,. \cr \endsplit$$
Using the second inequality with $y=0$, we get the first inequality of the theorem. 
Using the first inequality above if $y \in [1/2,1]$ and the second inequality above if $y \in [0,1/2]$ we get 
$$|f^{\prime}(y)|^2 \leq 27^2 (1 +\varepsilon_n)^2 \, n^3
\left( \sum_{k=1}^n{\left( \left( \frac{\lambda_k}{9n} \right)^2 + 4(k-1)^2\right)} \right)
\frac{1}{9n} \, \|f\|_{[0,1]}^2,$$
and the first statement of the theorem follows. 
\qed \enddemo

\demo{Proof of Theorem 9.2}
Let  $Q_n \in {\Cal P}_n$ defined by $Q_n(x) = T_n(2x-1)$, where
$T_n$ is the Chebyshev polynomial of degree $n$ on $[-1,1]$
defined by $T_n(\cos \theta) = \cos(n\theta)\,.$ As
$$|P_n^{\prime}(0)| = 2n^2 = 2n^2 \|P_n\|_{[0,1]}\,,$$
the theorem follows from (13.4).
\qed \enddemo

\demo{Proof of Theorem 10.1}
This follows from Lemma 12.12 by the substitution $x = e^{-t}$.
\qed \enddemo

\demo{Proof of Theorem 10.2}
This follows from Lemma 10.1 immediately.
\qed \enddemo

\demo{Proof of Theorem 11.1}
Observe that if $0 \neq f \in {\Cal E}_n^-$ is of the form
$$f(t) = \sum_{j=1}^n{a_j e^{\lambda_jt}}, \quad a_j, \lambda_j \in {\Bbb C}\,, \enskip \text {\rm Re}(\lambda_j) < 0\,,$$
then $g \in {\Cal E}_n$ defined by $g(t) = f(t)e^{t/2}$ is of the form 
$$g(t) = \sum_{j=1}^n{a_j e^{\lambda_jt}}, \quad a_j, \lambda_j \in {\Bbb C}\,, \enskip \text {\rm Re}(\lambda_j) < 1/2\,.$$
Now an application of Theorem 10.1 to $g$ gives
$$\frac{\left\|(\left(f^\prime(t)e^{t/2} + \frac 12 \, f(t)e^{t/2})\right)e^{-t/2} \right\|_{L_2[0,\infty)}}
{\|f(t)e^{t/2}e^{-t/2}\|_{L_2[0,\infty)}} = 
\frac{\|(g^\prime(t)e^{-t/2}\|_{L_2[0,\infty)}}
{\|g(t)e^{-t/2}\|_{L_2[0,\infty)}} $$
$$\leq \max_{1 \leq j \leq n}{\left|\lambda_j + \frac 12 \right|} +
\left( \sum_{j=1}^n{\left(1 - 2{\text {\rm Re}}\left(\lambda_j + \frac 12 \right)\right) 
\sum_{k=j+1}^n{\left(1 - 2{\text {\rm Re}}\left(\lambda_k + \frac 12\right) \right)}}\right)^{1/2} \,,$$	 
hence 
$$\frac{\|f^\prime\|_{L_2[0,\infty)}}{\|f\|_{L_2[0,\infty)}}
\leq \frac 12 + \max_{1 \leq j \leq n}{\left|\lambda_j + \frac 12 \right|} +
2 \left( \sum_{j=1}^n{{\text {\rm Re}}(\lambda_j)
\sum_{k=j+1}^n{{\text {\rm Re}}(\lambda_k)}}\right)^{1/2} \,.$$
\qed \enddemo

\medskip

\head 14. Appendix \endhead
The paper is self-contained without the results listed in this section. 
The results below are closely related to our new results in this paper. 
Theorems 14.1--14.6 have been proved by subtle Descartes system methods which 
can be employed in the case of exponential sums with only real exponents but not 
in the case of complex exponents. The reader may find it useful to compare 
the results in this section with the new results of the paper.     

Associated with a set of $\Lambda_n := \{\lambda_0, \lambda_1, \ldots, \lambda_n\}$ of distinct
real numbers let
$$E(\Lambda_n) :=
\text{\rm span}\{e^{\lambda_0t}, e^{\lambda_1t}, \ldots, e^{\lambda_nt}\} = 
\left \{f: f(t) = \sum_{j=0}^n{a_je^{\lambda_jt}} \,, \enskip  a_j \in {\Bbb R} \right \}\,.$$
The following result was proved in [14].

\proclaim{Theorem 14.1}
Suppose $\Lambda_n := \{\lambda_0,\lambda_1, \ldots,\lambda_n\}$ is a set of distinct
nonnegative real numbers. Let $0 < q \leq p \leq \infty$. Let $\mu$ be a non-negative integer.
There are constants $c_1 = c_1(p,q,\mu) > 0$ and $c_2 = c_2(p,q,\mu)$ depending only on
$p$, $q$, and $\mu$ such that
$$c_1 \left(\sum_{j=0}^n{\lambda_j}\right)^{\mu + 1/q - 1/p} \leq
\sup_{0 \neq f \in E(\Lambda_n)}
{\frac{\|f^{(\mu)}\|_{L_p(-\infty,0]}}{\|f\|_{L_q(-\infty,0]}}}
\leq c_2 \left(\sum_{j=0}^n{\lambda_j}\right)^{\mu + 1/q - 1/p}\,,$$
where the lower bound holds for all $0 < q \leq p \leq \infty$ and $\mu \geq 0$,
while the upper bound holds when $\mu = 0$ and $0 < q \leq p \leq \infty$,
and when $\mu \geq 1$, $p \geq 1$, and $0 < q \leq p \leq \infty$.
Also, there are constants $c_1 = c_1(q,\mu) > 0$ and $c_2 = c_2(q,\mu)$ depending only on
$q$ and $\mu$ such that
$$c_1 \left(\sum_{j=0}^n{\lambda_j}\right)^{\mu + 1/q} \leq
\sup_{0 \neq f \in E(\Lambda_n)} {\frac{|f^{(\mu)}(y)|}{\|f\|_{L_q(-\infty,y]}}}
\leq c_2 \left(\sum_{j=0}^n{\lambda_j}\right)^{\mu + 1/q}$$
for all $0 < q \leq \infty$, $\mu \geq 1$, and  $y \in {\Bbb R}$.
\endproclaim

In [15] we proved the following couple of theorems.

\proclaim{Theorem 14.2}
Suppose $\Lambda_n := \{\lambda_0,\lambda_1, \ldots,\lambda_n\}$ is a set of distinct
real numbers. Let $0 < q \leq p \leq \infty$, $a,b \in {\Bbb R}$, and $a<b$.
There are constants $c_3 = c_3(p,q,a,b) > 0$ and $c_4 = c_4(p,q,a,b)$
depending only on $p$, $q$, $a$, and $b$ such that
$$c_3 \left(n^2 + \sum_{j=0}^n{|\lambda_j|}\right)^{1/q - 1/p} \leq
\sup_{0 \neq f \in E(\Lambda_n)} {\frac{\|f\|_{L_p[a,b]}}{\|f\|_{L_q[a,b]}}}
\leq c_4 \left(n^2 + \sum_{j=0}^n{|\lambda_j|}\right)^{1/q - 1/p}\,.$$
\endproclaim

\proclaim{Theorem 14.3}
Suppose $\Lambda_n := \{\lambda_0,\lambda_1, \ldots,\lambda_n\}$ is a set of distinct 
real numbers. Let $0 < q \leq p \leq \infty$, $a,b \in {\Bbb R}$, and $a<b$.
There are constants $c_5 = c_5(p,q,a,b) > 0$ and $c_6 = c_6(p,q,a,b)$ depending only on 
$p$, $q$, $a$, and $b$ such that
$$c_5 \left(n^2 + \sum_{j=0}^n{|\lambda_j|}\right)^{1 + 1/q - 1/p} 
\leq \sup_{0 \neq f \in E(\Lambda_n)}{\frac{\|f^\prime\|_{L_p[a,b]}}{\|f\|_{L_q[a,b]}}}
\leq c_6 \left(n^2 + \sum_{j=0}^n{|\lambda_j|}\right)^{1 + 1/q - 1/p}\,,$$
where the lower bound holds for all $0 < q \leq p \leq \infty$,
while the upper bound holds when $p \geq 1$ and $0 < q \leq p \leq \infty$.
\endproclaim

Using the $L_{\infty}$ norm on a fixed subinterval $[a + \delta, b - \delta] \subset [a,b]$
in the numerator in Theorem 14.2, we proved the following essentially sharp result in [6].
For the sake of brevity let
$$\|f\|_A := \sup_{t \in A}{|f(t)|}$$
for a complex-valued function $f$ defined on a set $A \subset {\Bbb R}$.

\proclaim{Theorem 14.4} If $\Lambda_n := \{\lambda_0, \lambda_1, \ldots, \lambda_n\}$
is a set of distinct real numbers, then the inequality
$$\|f\|_{[a + \delta, b- \delta]} \leq e8^{1/p} \left(\frac{n+1}{\delta}\right)^{1/p}\|f\|_{L_p[a,b]}$$
holds for every $f \in E(\Lambda_n)$, $p > 0$, and $\delta \in \left(0, \textstyle{\frac 12}(b-a)\right)$.
\endproclaim

The key to this result is the following Remez-type inequality proved also in [6].
For the sake of brevity let
$$E_n: = \Biggl\{f : f(t) = a_0 + \sum^n_{j=1} {a_j e^{\lambda_j t}}\,, \enskip a_j, \lambda_j \in {\Bbb R} \Biggr\}$$
and
$$E_n(s) := \{f \in E_n: m \left( \left\{x \in [-1,1]: |f(x)| \leq 1 \right\} \right) \geq 2-s\}\,,$$
where $m(A)$ denotes the Lebesgue measure of a measurable set $A \subset {\Bbb R}$.

\proclaim{Theorem 14.5}
Let $s \in \left( 0, \frac 12 \right]\,.$
There are absolute constants $c_7 > 0$ and $c_8 > 0$ such that
$$\exp(c_7\min\{ns,(ns)^2\}) \leq \sup_{f \in E_n(s)}{|f(0)|} \leq \exp(c_8\min\{ns, (ns)^2\})\,.$$
\endproclaim

An essentially sharp Bernstein-type inequality for $E_n$ is proved in [4].   

\proclaim{Theorem 14.6}
We have
$$\frac {1}{e-1}\,\frac {n-1}{\min\{y-a,b-y\}} \leq \sup_{0 \neq f \in E_n} \frac{|f^\prime(y)|}{\|f\|_{[a,b]}}
\leq \frac {2n - 1}{\min\{y-a,b-y\}}\,, \qquad y\in(a,b)\,.$$
\endproclaim

Having real exponents $\lambda_j$ in Theorems 1.1--1.6 is essential in the proofs using subtle Descartes 
system methods. There are other important inequalities proved for the classes $E(\Lambda_n)$ associated with a set 
$\Lambda_n := \{\lambda_0, \lambda_1, \ldots, \lambda_n\}$ of distinct real exponents. See [5], 
for instance, where the proofs are using Descartes system methods as well.   

Let $V_n$ be a vector space of complex-valued functions defined on ${\Bbb R}$ of
dimension $n+1$ over ${\Bbb C}$. We say that $V_n$ is shift invariant (on ${\Bbb R}$)
if $f \in V_n$ implies that $f_a \in V_n$ for every $a \in {\Bbb R}$, where
$f_a(x) := f(x-a)$ on ${\Bbb R}$.
Associated with a set of $\Lambda_n := \{\lambda_0, \lambda_1, \ldots, \lambda_n\}$ of distinct
COMPLEX numbers let 
$$E^c(\Lambda_n) :=
\text{\rm span}\{e^{\lambda_0t}, e^{\lambda_1t}, \ldots, e^{\lambda_nt}\} =
\left \{f: f(t) = \sum_{j=0}^n{a_je^{\lambda_jt}}, \, \enskip a_j \in {\Bbb C} \right \}\,.$$
Elements of $E^c(\Lambda_n)$ are called exponential sums of $n+1$ terms.
Examples of shift invariant spaces of dimension $n+1$ include $E^c(\Lambda_n)$.
In [7] we proved a result analogous to Theorem 14.4 for complex exponents $\lambda_j$, 
in which case Descartes system methods cannot help us in the proof.

\proclaim{Theorem 14.7} Let $V_n \subset C[a,b]$ be a shift invariant vector space of complex-valued
functions defined on ${\Bbb R}$ of dimension $n+1$ over ${\Bbb C}$. Let $p \in (0,2]$. Then
$$\|f\|_{[a + \delta,b-\delta]} \leq
2^{2/p^2} \left(\frac{n+1}{\delta} \right)^{1/p} \|f\|_{L_p[a,b]}$$
for every $f \in V_n$, $p \in (0,2]$, and $\delta \in \left (0,\frac 12(b-a) \right)\,,$ and
$$\|f\|_{[a + \delta,b-\delta]} \leq 2^{1/2} \left( \frac{n+1}{\delta} \right)^{1/2} (b-a)^{(p-2)/p}\|f\|_{L_p[a,b]}$$ 
for every $f \in V_n$, $p \geq 2$, and $\delta \in \left (0,\frac 12(b-a) \right)\,.$
\endproclaim

It is well known by considering the the case of algebraic polynomials of degree $n$ 
that, in general, the size of the factor $(n+1)^{1/p}$ in Theorem 14.7 cannot be improved for $p \in (0,2]$. 
On the other hand for $p \geq 2$ the size of the factor $(n+1)^{1/2}$ in 
the inequality 
$$\split \|f\|_{[a + \delta,b-\delta]} \leq & \, 2^{1/2} \left(\frac{n+1}{\delta} \right)^{1/2} \|f\|_{L_2[a,b]} \cr
\leq & \, 2^{1/2} \left(\frac{n+1}{\delta} \right)^{1/2} (b-a)^{(p-2)/(2p)}\|f\|_{L_p[a,b]} \cr \endsplit$$
cannot be improved. This can be seen by taking lacunary trigonometric polynomials. See the theorem below from 
[30, p. 215].

\proclaim{Theorem 14.8}
Let $(k_j)$ be a strictly increasing sequence of nonnegative integers satisfying
$$k_{j+1} > \alpha k_j\,, \qquad j=1,2,\ldots\,,$$
where $\alpha > 1$. Let 
$$Q_n(t) = \sum_{j=1}^n{\cos(2\pi k_j(t-\theta_j))}\,.$$
Then for every $q > 0$ there are constants $A_{q,\alpha} > 0$ and $B_{q,\alpha} > 0$ depending only on 
$q$ and $\alpha$ such that 
$$A_{q,\alpha} n^{1/2} \leq \|Q_n\|_{L_q[0,2\pi]} \leq B_{q,\alpha} n^{1/2}$$
for every $n \in {\Bbb N}$ and $q > 0$.
\endproclaim

\head 15. Acknowledgements \endhead
The author wishes to thank Sergey Denisov for the motivating questions and e-mail 
discussions related to the paper. The author wishes to thank Stephen Choi as well 
for reading earlier versions of my paper carefully, pointing out many misprints, 
and his advise about the structure of the presentation.

\enddocument